
%

\input amstex.tex
\loadmsam
\loadmsbm
\loadbold
\input amssym.tex
\baselineskip=13pt plus 2pt
\documentstyle{amsppt}
\pageheight{45pc}
\pagewidth{33pc}

\magnification=1200
\overfullrule=0pt

\def\enddemos{\hfill$\square$\enddemo}
\def\hb{\hfil\break}
\def\n{\noindent}
\def\ov{\overline}
\def\smatrix{\smallmatrix}

\def\pmatrix{\left(\smatrix}
\def\endpmatrix{\endsmallmatrix\right)}
\def\upi{\pmb{\pi}}
\def\pii{\pmb{\pi}}
\def\upsi{\pmb{\psi}}

\def\A{\Bbb A}
\def\C{\Bbb C}
\def\F{\Bbb F}
\def\G{\Bbb G}

\def\Z{\Bbb Z}

\def\1{\bold 1}
\def\r{\bold r}
\def\s{\bold s}
\def\u{\bold u}
\def\v{\bold v}
\def\w{\bold w}

\def\uC{\bold C}
\def\uG{\bold G}
\def\uT{\bold T}
\def\Gal{\operatorname{Gal}}

\def\max{\operatorname{max}}
\def\Vol{\operatorname{vol}}
\def\vol{\operatorname{vol}}
\def\GL{\operatorname{GL}}

\def\GSp{\operatorname{GSp}}
\def\SL{\operatorname{SL}}
\def\PGL{\operatorname{PGL}}

\def\SO{\operatorname{SO}}

\def\Ind{\operatorname{Ind}}
\def\tr{\operatorname{tr}}
\def\Re{\operatorname{Re}}
\def\Lie{\operatorname{Lie}}
\def\Ad{\operatorname{Ad}}

\def\diag{\operatorname{diag}}
\def\fz4{\text{\rm ([FZ4])}}

\leftheadtext{Yuval Z. Flicker and Dmitrii Zinoviev}
\rightheadtext{Twisted character of a small representation}
\topmatter
\title Twisted character of a small representation of GL(4)\endtitle
\author Yuval Z. Flicker and Dmitrii Zinoviev\endauthor
\footnote"~"{\n Department of Mathematics, The Ohio State University,
231 W. 18th Ave., Columbus, OH 43210-1174;\hb
email: flicker\@math.ohio-state.edu. Partially supported by a Lady Davis
Visiting Professorship at the Hebrew University, Max-Planck scholarship
at MPI Bonn, and the Humboldt Stiftung.\hb
\indent Institute for Problems in Information Transmission, Russian Academy
of Sciences, Bolshoi Karetnyi per. 19, GSP-4, Moscow 101447, Russia; 
email: zinov\@iitp.ru. Partially supported by Russian Academy of Science
research grant 03-01-00098.\hb
}

\abstract We compute by a purely local method the (elliptic)
$\theta$-twisted character $\chi_{\pi_Y}$ of the representation
$\pi_Y=I_{(3,1)}(1_3\times\chi_Y)$ of $G=\GL(4,F)$, where $F$ is a
$p$-adic field, $p\not=2$, and $Y$ is an unramified quadratic
extension of $F$; $\chi_Y$ is the nontrivial character of
$F^\times/N_{Y/F}Y^\times$.

The representation $\pi_Y$ is normalizedly induced from
$\pmatrix m_3&\ast\\ 0&m_1\endpmatrix\mapsto
\chi_Y(m_1)$, $m_i\in\GL(i,F)$, on the maximal parabolic
subgroup of type $(3,1)$; $\theta$ is the ``transpose-inverse''
involution of $G$.

We show that the twisted character $\chi_{\pi_Y}$ of $\pi_Y$ is an
unstable function: its value at a twisted regular elliptic conjugacy
class with norm in $C_Y=\uC_Y(F)$=``$(\GL(2,Y)/F^\times)_F$'' is minus
its value at the other class within the twisted stable conjugacy class.
It is $0$ at the classes without norm in $C_Y$. Moreover $\pi_Y$ is the
endoscopic lift of the trivial representation of $C_Y$.

We deal only with unramified $Y/F$, as globally this case occurs almost
everywhere. The case of ramified $Y/F$ would require another paper.

Our $\uC_Y=$``$(\text{R}_{Y/F}\GL(2)/\GL(1))_F$'' has $Y$-points $\uC_Y(Y)=
\{(g,g')\in \GL(2,Y)\times\GL(2,Y);$ $\det(g)=\det(g')\}/Y^\times$
($Y^\times$ embeds diagonally); $\sigma(\not= 1)$ in $\Gal(Y/F)$
acts by $\sigma(g,g')=(\sigma g',\sigma g)$. It is a $\theta$-twisted
elliptic endoscopic group of GL(4).

Naturally this computation plays a role in the theory of lifting of $C_Y$
and GSp(2) to GL(4) using the trace formula, to be discussed elsewhere.

Our work extends -- to the context of nontrivial central characters
-- the work of [FZ4], where representations of $\PGL(4,F)$ are studied.
In [FZ4] a 4-dimensional analogue of the model of the small
representation of $\PGL(3,F)$ introduced with Kazhdan in [FK] in
a 3-dimensional case is developed, and the local method of computation
introduced in [FZ3] is extended. As in [FZ4] we use here the
classification of twisted (stable) regular conjugacy classes in
$\GL(4,F)$ of [F], motivated by Weissauer [W].

\endabstract
\endtopmatter

\document
\heading Introduction\endheading
Let $\pi$ be an admissible representation (see Bernstein-Zelevinsky [BZ],
2.1) of a $p$-adic reductive group $G$. Its character $\chi_\pi$ is a
complex valued function defined by $\tr\pi(fdg)=\int_G\chi_\pi(g)f(g)dg$
for all complex valued smooth compactly supported measures $fdg$ ([BZ],
2.17). It is smooth on the regular set of the group $G$. The character
is important since it characterizes the representation up to equivalence.
A fundamental result of Harish-Chandra [H] establishes that the character
is a locally integrable function in characteristic zero.

Let $\theta$ be an automorphism of finite order of the group $G$. Define
${}^\theta\pi$ by ${}^\theta\pi(g)=\pi(\theta(g))$. When $\pi$ is invariant
under the action of $\theta$ (thus ${}^\theta\pi$ is equivalent to $\pi$),
Shintani and others introduced an extension of $\pi$ to the semidirect
product $G\rtimes\langle\theta\rangle$. The twisted character
$\chi_\pi(g\times\theta)$ is defined by $\tr\pi(fdg\times\theta)
=\int_G\chi_\pi(g\times\theta)f(g)dg$ for all $fdg$. It depends only
on the $\theta$-conjugacy class $\{hg\theta(h)^{-1};h\in G\}$ of $g$.
It is again smooth on the $\theta$-regular set, and characterizes the
$\theta$-invariant irreducible $\pi$ up to isomorphism. Moreover, it
is locally integrable (see Clozel [C]) in characteristic zero.



Characters provide a very precise tool to express a relation of
representations of different groups, called lifting, initiated by
Shintani and studied extensively in the case of base change, and
also in non base change situations such as twisting by characters
(Kazhdan [K], Waldspurger [Wa]), and the symmetric square lifting
from SL(2) to PGL(3) ([Fsym], [FK]). In this last case twisted
characters of $\theta$-invariant representations of PGL(3) are
related to packets of representations of SL(2), and $\theta$ is
the involution sending $g$ to its transpose-inverse.

{\it The aim of the present work is to compute the twisted $($by $\theta)$
character of a specific representation $\pi=\pi_Y=I_{(3,1)}(1_3\times\chi_Y)$,
of the group $G=\GL(4,F)$, $F$ a $p$-adic field, $p$ odd}. Here $Y/F$ is an
unramified quadratic extension and $\chi_Y$ is the quadratic character of
$F^\times$ which is trivial on the group $N_{Y/F}Y^\times$, where $N_{Y/F}$
is the norm map from $Y$ to $F$. This $\pi$ is normalizedly induced from the
representation $\pmatrix m_3&\ast\\ 0&m_1\endpmatrix\mapsto\chi_Y(m_1)$,
$m_i\in\GL(i,F)$, of the standard (upper triangular) maximal parabolic
subgroup $P$ of type $(3,1)$. It is invariant under the involution
$\theta(g)=J^{-1}{}^tg^{-1}J$, where
$J=(a_i\delta_{i,5-j})$, $a_1=a_2=1$, $a_3=a_4=-1$.

A natural setting for the statement of our result is the theory of
liftings to the group $\uG=\GL(4)$ from its $\theta$-twisted
endoscopic (see Kottwitz-Shelstad [KS]) $F$-group
$$\uC_Y=\{(g,g')\in\GL(2)\times\GL(2);\,\det g=\det g'\}/\G_m,$$
where the multiplicative group $\G_m=\GL(1)$
embeds as $z\mapsto(zI_2,zI_2)$, $I_2$ is the identity $2\times 2$
matrix, with $\Gal(\ov{F}/F)$-action which is a composition of the usual
Galois action on each of the two factors GL(2) with the transposition
$(g,g')\mapsto(g',g)$ if $\sigma\in\Gal(\ov{F}/F)$ has nontrivial
restriction to $Y$. Here $\ov{F}$ is a separable algebraic closure
of $F$ containing $Y$.

The corresponding map $\lambda_Y$ of dual groups is simply the natural
embedding in $\hat{G}=\GL(4,\C)$ of the non connected
$\hat{C}_Y=Z_{\hat{G}}(\hat{s}\hat\theta)$
$$=\left\{g\in\hat{G}=\GL(4,\C);\,g\hat{s}J{}^tg=\hat{s}J=
\pmatrix 0 &\omega\\ \omega^{-1} & 0\endpmatrix\right\}=\operatorname{O}
\left(\pmatrix 0 &\omega\\ \omega^{-1} & 0\endpmatrix,\C\right)$$
$$=\left\langle\pmatrix aB & bB\\ cB & dB\endpmatrix,\,
\iota=\pmatrix 1&&&0\\&0&1&\\&1&0&\\0&&&1\endpmatrix;\,
\left(A=\pmatrix a & b\\ c & d\endpmatrix,B\right)
\in\left\{\GL(2,\C)^2;\,\det A\cdot\det B=1\right\}/\C^\times\right\rangle,$$
where $z\in\C^\times$ embeds as the central element $(z,z^{-1})$, and
where $\hat{s}=\diag(-1,1,-1,1)$ and $\omega=\pmatrix 0&-1\\ 1&0\endpmatrix$.
Thus $\hat{C}_Y$ is the $\hat\theta$-centralizer in $\hat{G}$ of the
semisimple element $\hat{s}$ (and $\hat\theta$ is defined on $\hat{G}$
by the same formula that defines $\theta$ on $G$), and $\Gal(Y/F)$ acts
via conjugation by $\iota$.

Indeed, our result can be viewed as asserting that the $\theta$-invariant
representation $\pi$ of $G=\uG(F)$, whose central character is
$\chi_Y\not=1$ of order two, is the endoscopic lift of the trivial
representation of $C_Y=\uC_Y(F)=(\GL(2,Y)/F^\times)_F$. The subscript
$F$ here indicates that $(\GL(2,Y)/F^\times)_F$ consists of $g$ in
$\GL(2,Y)/F^\times$ with $\det(g)$ in $F^\times/F^{\times 2}$.

To state this we note that the embedding $\lambda_Y:\hat{C}_Y\to\hat{G}$
defines a norm map. This norm map relates the stable $\theta$-conjugacy
classes in $G$ with stable conjugacy classes in $C_Y$, where ``stable''
means the elements in $G$ of an orbit in the points of $\uG$ in a separable
algebraic closure $\ov{F}$ of $F$. The crucial case is that of
$\theta$-elliptic elements. A stable $\theta$-conjugacy class
consists of several $\theta$-conjugacy classes. The stable $\theta$-conjugacy
classes of elements in $G$, and the $\theta$-conjugacy classes within the
stable $\theta$-classes, have been described recently in [F], in analogy
with the description of the (stable) conjugacy classes in the group of
symplectic similitudes $\GSp(2,F)$ of Weissauer [W]. In fact in [FZ4] the
simpler case of $\PGL(4,F)$ is used, but here, as in [F], we deal
with $\theta$-classes in $\GL(4,F)$. We give here full details of the
description in our case.

There are four types of $\theta$-elliptic elements of $G$, named in [F]
and here I, II, III, IV, depending on their splitting behaviour. As in
[F], our work relies on an explicit presentation of representatives of
the $\theta$-conjugacy classes within the stable such classes in $G$.
We present here the same set of representatives as in [FZ4].

The norm map, which we describe explicitly here, relates $\theta$-conjugacy
classes of types II and IV to conjugacy classes in $C_Y$. It does not relate
classes of types I, III to classes in $C_Y$. Our ``quadratic'' case behaves
then in a complementary fashion to that of [FZ4], where $\theta$-conjugacy
classes of types I, III are related to conjugacy classes in the group
$C=\SO(4)$ of [FZ4], but $\theta$-conjugacy classes of types II, IV are not
related to conjugacy classes in $C$.


The stable $\theta$-conjugacy classes of types II and IV come associated
with a quadratic extension $E/E_3$, where $Y=E_3$ is a quadratic extension
of $F$. The two $\theta$-conjugacy classes $g_r$ within the stable
$\theta$-classes are parametrized by $r$ in $E_3^\times/N_{E/E_3}E^\times$.
We prove

\proclaim{Theorem} The value of the $\theta$-character
$\chi_\pi(g\times\theta)$ of $\pi=\pi_Y$ at the $\theta$-regular element
$g=g_r$ of type II or IV, multiplied by a suitable Jacobian
${{\Delta(g_r\theta)}\over{\Delta_C(Ng)}}$, is $2\kappa(r);$ here $\kappa$
is the $\not=1$ character of $E_3^\times/N_{E/E_3}E^\times$.
At any $\theta$-regular element $g$ of type I or III,
$\chi_\pi(g\times\theta)=0$.\endproclaim

In particular the character $\chi_\pi(g\times\theta)$ is
an unstable function, namely its value at one $\theta$-conjugacy class
within a stable $\theta$-conjugacy class of type II or IV is negative its
value at the other $\theta$-conjugacy class.

We deal only with unramified $Y/F$, as globally this case occurs almost
everywhere. The case of ramified $Y/F$ would require another paper.

Our result is a special case of the lifting with respect to $\lambda_Y$
to the group $G=\GL(4,F)$ of representations of the group $C_Y=(\GL(2,Y)
/F^\times)_F$.

Our work develops the method of [FZ4] to the context of representations
with nontrivial central characters. We use a model of our representation
$\pi_Y=I_{(3,1)}(1_3\times\chi_Y)$, different from the standard model of
a parabolically induced representation. It is a twist of the four dimensional
analogue of [FZ4], dealing with $\pi_4=I_{(3,1)}(1_3\times 1)$, of a (three
dimensional) model introduced and used with Kazhdan in [FK] to compute
the twisted (by transpose-inverse) character of the representation
$\pi_3=I_{(2,1)}(1_2)$ of $\PGL(3,F)$ normalizedly induced from the trivial
representation of the maximal parabolic subgroup. We do not use our results
to prove the fundamental lemma since in our case, as well as that of [FZ4],
the fundamental lemma is already established in [F]. In the case of the
symmetric square lifting from $\SL(2,F)$ to $\PGL(3,F)$, an analogous
purely local and simple proof of the fundamental lemma was given in [Fsym.
Unit elements].

The work of [FK] uses local arguments to compute the twisted character of
$\pi_3$ on one of the two twisted conjugacy classes within the stable one
(where the quadratic form is anisotropic), and global arguments to reduce
the computation on the other class (where the quadratic form is isotropic)
to that computed by local means. A purely local computation for the second
class is given in [FZ3]. In [FZ4] this local computation is developed in a
four dimensional projective case. A global type of argument as in [FK] is
harder to apply as there are not enough anisotropic quadratic forms in the
four dimensional case. Anyway, a simpler, local proof, is better. Here we
extend the work of [FZ4] to $\theta$-invariant representations of $\GL(4,F)$
whose central character is nontrivial, necessarily quadratic. Our work is
parallel to -- but entirely independent of -- the work of [FZ4].

\heading Conjugacy classes\endheading

Let $F$ be a local nonarchimedean field, and $R$ its ring of integers.
Put $\uG=\GL(4)$, $G=\uG(F)$ and $K=\uG(R)$. Put
$\uC_Y=\{(g_1,g_2)\in \GL(2)\times \GL(2);\,\det(g_1)=\det(g_2)\}/\G_m$
($\G_m$ embeds diagonally), viewed as a group over $F$ with Galois action
$\tau(g,g')=(\tau g,\tau g')$ unless $\tau\in\Gal(\ov{F}/F)$ has nontrivial
restriction to $Y$, in which case $\tau(g,g')=(\tau g',\tau g)$, where
$\tau(g_{ij})=(\tau g_{ij})$. It is a form of the group $\uC$ of [FZ4],
and in particular $\uC_Y(Y)=\uC(Y)$, but the $\Gal(\ov{F}/F)$-action is
different: $\tau\in\Gal(\ov{F}/F)$ takes $(g,g')$ of $\uC(\ov{F})$ to
$(\tau g,\tau g')$. Then $C_Y=\uC_Y(F)=\{g\in GL(2,Y)/F^\times;\,
\det(g)\in F^\times\}$ and $K_{C_Y}=\uC_Y(R)$.
Set $\theta(\delta)=J^{-1}{}^t\delta^{-1}J$ for $\delta$ in $G$.
Here $J$ is $\pmatrix 0 & w \\ -w & 0 \endpmatrix$, where
$w\,=\,\pmatrix 0 & 1 \\ 1 & 0 \endpmatrix$.
Fix a separable algebraic closure $\ov{F}$ of $F$.
The elements $\delta$, $\delta'$ of $G$ are called ({\it stably})
$\theta$-{\it conjugate} if there is $g$ in $G$ (resp. $\GL(4,\ov{F}))$
with $\delta'=g^{-1}\delta\theta(g)$.

Results of [F] concerning (stable) $\theta$-twisted regular conjugacy 
classes are recalled in [FZ4], pp. 337-338. There are four types of 
$\theta$-elliptic classes, but the norm map $N$ from $G$ to $C_Y$ 
relates only the twisted classes in $G$ of type II and IV to conjugacy 
classes in $C_Y$. We should then expect the twisted character of the 
representation considered here to vanish on the twisted classes of 
type I and III. 

\heading Norm map\endheading

The norm map $N:\uG\to\uC_Y$ is defined on the diagonal torus $\uT^\ast$
of $\uG$ by $$N(\diag(a,b,c,d))=(\diag(ab,cd),\diag(ac,bd)).$$ Since
both components have determinant $abcd$, the image of $N$
is indeed in $\uC_Y$.

In type II we have $a\in E_1^\times$, $E_1=E^\tau=F(\sqrt D)$,  $b\in
E_2^\times$, $E_2=E^{\sigma\tau}=F(\sqrt{AD})$, and the norm map becomes
$$N(\diag(a,b,\tau b,\sigma a))=(\diag(ab,\tau(b)\sigma(a)),
\diag(a\tau(b),b\sigma(a))).$$
The two components on the right are mapped to each other by $\tau$, while
the pairs of eigenvalues ($\{ab,\tau(b)\sigma(a)\}$ and $\{a\tau(b),
b\sigma(a)\}$) are permuted by $\sigma$. Hence the right side defines a
conjugacy class in $\GL(2,E_3)_F$ (the determinant $ab\cdot\tau(b)\sigma(a)
=a\tau(b)\cdot b\sigma(a)$ lies in $F^\times$, and $E_3$ is the fixed field
of $\sigma$ in $E$). We choose $Y$ to be the quadratic extension $E_3$ of
$F$. The image of this torus is the torus (up to conjugacy) in $C_Y(F)$
which splits over the biquadratic extension $E$ of $F$.

In type IV we have $\alpha\in E^\times$, $E=E_3(\sqrt D)$,
$\alpha\sigma^2\alpha\in E_3^\times$, $E_3=F(\sqrt A)$,
and the norm map becomes
$$N(\diag(\alpha,\sigma\alpha,\sigma^3\alpha,\sigma^2\alpha))=
(\diag(\alpha\sigma\alpha,\sigma^2\alpha\sigma^3\alpha),
\diag(\alpha\sigma^3\alpha,\sigma\alpha\sigma^2\alpha)).$$
Here $\sigma^3$ permutes the two diagonal matrices on the right,
and $\sigma^2$ permutes each pair of eigenvalues.
Since both components of $N(\ast)$ have equal determinants in $F^\times$,
$\diag(\ast,\ast)$ defines a conjugacy class in $\GL(2,E_3)_F$. Hence
the norm map defines a conjugacy class in $C_Y=\uC_Y(F)$ for each
$\theta$-stable conjugacy class of type IV in $G=\uG(F)$, where we
take $Y$ to be $E_3$.

In types I and III the image of the map $N$ does not correspond to any
conjugacy class in $C_Y$, for any quadratic extension $Y$ of $F$.

\heading Jacobians\endheading

The character relation that we study relates the product of the value
at $t$ of the twisted character of our representation $\pi_Y=I_{(3,1)}(1_3
\times\chi_Y)$ by a factor $\Delta(t\times\theta)$, with the product by a
factor $\Delta_C(Nt)$ of the value at $Nt$ of the character of the (trivial)
representation $\1_{C_Y}$ of $C_Y$ which lifts to $\pi_Y$.

The factors $\Delta(t\times\theta)$ and $\Delta_C(Nt)$ are defined 
and computed in [FZ4], pp. 339-340. We have
$${\Delta(t\theta)\over\Delta_C(Nt)}=\left|{(a-d)^2\over{ad}}\cdot
{(b-c)^2\over{bc}}\right|^{1/2}.$$

Then in case II if $t=\diag(a,b,\tau b,\sigma a)$, $a=a_1+a_2\sqrt D\in
E_1^\times$, $b=b_1+b_2\sqrt{AD}\in E_2^\times$, we get
$${\Delta(t\theta)\over\Delta_C(Nt)}
=\left|{(a-\sigma a)^2\over{a\sigma a}}\cdot
{(b-\sigma b)^2\over{b\sigma b}}\right|^{1/2}
=\left|{(2a_2\sqrt D)^2\over{a_1^2-a_2^2D}}
\cdot{(2b_2\sqrt{AD})^2\over{b_1^2-b_2^2AD}}\right|^{1/2}.$$

In case IV, if $t=\diag(\alpha,\sigma\alpha,\sigma^3\alpha,\sigma^2\alpha)$,
$\alpha=a+b\sqrt D$, $a=a_1+a_2\sqrt A$, $b=b_1+b_2\sqrt A$,
$\sigma\alpha=\sigma a+\sigma b\sqrt{\sigma D}$, $\sigma^3\alpha=\sigma a
-\sigma b\sqrt{\sigma D}$, $\alpha-\sigma^2\alpha=2b\sqrt D$,
$\sigma(\alpha-\sigma^2\alpha)=2\sigma b\sqrt{\sigma D}$, and
$${\Delta(t\theta)\over\Delta_C(Nt)}
=\left|{(\alpha-\sigma^2\alpha)^2\over{\alpha\sigma^2\alpha}}\cdot
{\sigma(\alpha-\sigma^2\alpha)^2\over{\sigma\alpha\sigma^3\alpha}}\right|^{1/2}
=\left|{(4b\sigma b)^2D\sigma D\over{(a^2-b^2D)\sigma(a^2-b^2D)}}\right|^{1/2}.
$$

\heading Characters\endheading

Denote by $f$ (resp. $f_{C_Y}$) a complex-valued
compactly-supported smooth (thus locally-constant
since $F$ is nonarchimedean) function on $G$
(resp. $C_Y$). Fix Haar measures on $G$ and on $C_Y$.

By a $G$-module $\pi$ (resp. $C_Y$-module $\pi_{C_Y}$)
we mean an admissible representation ([BZ]) of $G$ (resp.
$C_Y$) in a complex space. An irreducible $G$-module
$\pi$ is called {\it{$\theta$-invariant}} if it is
equivalent to the $G$-module $^\theta\pi$, defined by
$^\theta\pi(g)=\pi(\theta (g))$. In this case there is
an intertwining operator $A$ on the space of $\pi$ with
$\pi(g)A=A\pi(\theta (g))$ for all $g$. Since $\theta^2=1$
we have $\pi(g)A^2=A^2\pi(g)$ for all $g$, and since
$\pi$ is irreducible $A^2$ is a scalar by Schur's lemma.
We choose $A$ with $A^2=1$. This determines $A$ up to a
sign. When $\pi$ has a Whittaker model, which happens for
all components of cuspidal automorphic representations of
the adele group $\GL(4,\A)$, we specify a normalization of
$A$ which is compatible with a global normalization, as follows,
and then put $\pi(g\times\theta)=\pi(g)\times A$.

Fix a nontrivial character $\upsi$ of $F$ in $\C^\times$,
and a character $\psi(u)=\upsi(a_{1,2}+a_{2,3}-a_{3,4})$ of
$u=(u_{i,j})$ in the upper triangular subgroup $U$ of $G$.
Note that $\psi(\theta(u))=\psi(u)$. Assume that $\pi$ is a
{\it non degenerate} $G$-module, namely it embeds in the space
of ``Whittaker'' functions $W$ on $G$, which satisfy -- by
definition -- $W(ugk)=\psi(u)W(g)$ for all $g\in G$, $u\in U$,
$k$ in a compact open subgroup $K_W$ of $K$, as a $G$-module
under right shifts: $(\pi(g)W)(h)=W(hg)$. Then ${}^\theta\pi$
is non degenerate and can be realized in the space of functions
${}^\theta W(g)=W(\theta(g))$, $W$ in the space of $\pi$. We take
$A$ to be the operator on the space of $\pi$ which maps $W$ to
${}^\theta W$.

A $G$-module $\pi$ is called {\it unramified} if the
space of $\pi$ contains a nonzero $K$-fixed vector.
The dimension of the space of $K$-fixed vectors is
bounded by one if $\pi$ is irreducible. If $\pi$
is $\theta$-invariant and unramified, and $v_0\neq 0$
is a $K$-fixed vector in the space of $\pi$, then
$Av_0$ is a multiple of $v_0$ (since $\theta K=K$),
namely $Av_0=cv_0$, with $c=\pm 1$. Replace $A$ by
$cA$ to have $Av_0=v_0$, and put $\pi(\theta)=A$.

When $\pi$ is (irreducible) unramified and has a Whittaker
model, both normalizations of the intertwining operator are
equal. In this case $\upsi$ is unramified (trivial on $R$
but not on $\upi^{-1}R$, where $\upi$ is a generator of the
maximal ideal of $R$), and there exists a unique Whittaker
function $W_0$ in the space of $\pi$ with respect to $\psi$
with $W_0=1$ on $K$. It is mapped by $\pi(\theta)=A$ to ${}^\theta W_0$,
which satisfies ${}^\theta W_0(k)=1$ for all $k$ in $K$ since
$K$ is $\theta$-invariant. Namely $A$ maps the unique normalized
(by $W_0(K)=1$) $K$-fixed vector $W_0$ in the space of $\pi$ to the
unique normalized $K$-fixed vector ${}^\theta W_0$ in the space of
${}^\theta\pi$, and we have ${}^\theta W_0=W_0$.

For any (admissible) $\pi$ and (smooth) $f$ the convolution operator
$\pi(f dg)=\smallint_G f(g)\pi(g)dg$ has
finite rank. If $\pi$ is $\theta$-invariant put
$\pi(f dg\times\theta)=\smallint_G f(g)\pi(g)\pi(\theta)dg$.
Denote by $\tr~\pi(f dg\times\theta)$ the trace of
the operator $\pi(f dg\times\theta)$. It depends on
the choice of the Haar measure $dg$, but the
({\it{twisted}}) {\it{character}} $\chi_\pi$ of
$\pi$ does not; $\chi_\pi$ is a locally-integrable
complex-valued function on $G\times\theta$ (see [C], [H])
which is $\theta$-conjugacy invariant and
locally-constant on the $\theta$-regular set, with
$\tr\,\pi(f dg\times\theta)=\smallint_G
f(g)\chi_\pi(g\times\theta)dg$ for all $f$.

Local integrability is not used in this work; rather it is recovered
for our twisted character.

\heading Small representation\endheading

To describe the $G$-module of interest in this paper, take $P$ to be
the upper triangular parabolic subgroup of type (3,1), and fix its
Levi factor to be $M=\{m=\diag(m_3,m_1);\,m_3\in\GL(3,F),\,m_1\in F^\times\}$.
It is isomorphic to $\GL(3,F)\times F^\times$. Let $\delta$ denote (as above)
the character $\delta(p)=|\Ad(p)|\Lie N|$ of $P$; it is trivial on the
unipotent radical $N$ (= $F^3)$ of $P$. Then the value of $\delta$ at
$p=mn$ is $\vert m_1^{-3}\det m_3\vert$. Denote by $I(\pi_1)$ the $G$-module
$\pi=\Ind(\delta^{1/2}\pi_1; P,G)$ normalizedly induced from $\pi_1$ on
$P$ to $G$. It is clear from [BZ] that when $\pi_1$ is self-contragredient
and $I(\pi_1)$ is irreducible then $I(\pi_1)$ is $\theta$-invariant, and it
is unramified if and only if $\pi_1$ is unramified.

{\it Our aim in this work is to compute the $\theta$-twisted
character $\chi_{\pi_Y}$ of the $\GL(4,F)$-module
$\pi_Y=I_{(3,1)}(1_3\times\chi_Y)$, where $1_3\times\chi_Y$ is the
$P$-module $\pmatrix m_3&\ast\\ 0&m_1\endpmatrix\mapsto\chi_Y(m_1)$,
$\chi_Y$ is a quadratic character of $F^\times$, $m_i\in\GL(i,F)$,
by purely local means.}

We begin by describing a useful model of our representation, in analogy
with the models of [FK] and [FZ4] of analogous representations
$I_{(2,1)}(1_2)$ of $\PGL(3,F)$ and $I_{(3,1)}(1_3)$ of $\PGL(4,F)$.
Indeed we shall express $\pi_Y$ as an integral operator in a convenient
model, and integrate the kernel over the diagonal to compute the character
of $\pi_Y$.

Denote by $\mu=\mu_s$ the character $\mu(x)=\vert x\vert^{(s+1)/2}$ of
$F^\times$, and by $\chi_Y$ a quadratic character of $F^\times$. This
pair $(\mu,\chi_Y)$ defines a character $\mu_P=\mu_{s,Y,P}$ of $P$,
trivial on $N$, by $\mu_P(p)=\mu((\det m_3)/m_1{^3})\chi_Y(m_1)$
if $p=mn$ and $m=\pmatrix m_3&0\\0&m_1\endpmatrix$ with $m_3$ in
$\GL(3,F)$, $m_1$ in $\GL(1,F)$. If $s=0$, then $\mu_P=\delta^{1/2}\chi_Y$,
where viewed as a character on $P$, $\chi_Y$ takes the value $\chi_Y(m_1)$
at $p$. Let $W_s=W_s^Y$ be the space of complex-valued smooth functions
$\psi$ on $G$ with $\psi(pg)=\mu_P(p)\psi(g)$ for all $p$ in $P$ and $g$
in $G$. The group $G$ acts on $W_s$ by right translation:
$(\pi_s(g)\psi)(h)=\psi(hg)$. By definition, $I_{(3,1)}(1_3\times\chi_Y)$
is the $G$-module $W_s$ with $s=0$. The parameter $s$ is introduced for
purposes of analytic continuation.

We prefer to work in another model $V_s=V_s^Y$ of the
$G$-module $W_s$. Let $V$ denote the space of
column 4-vectors over $F$. Let $V_s$ be the space
of smooth complex-valued functions $\phi$ on
$V-\{0\}$ with $\phi(\lambda\v)=\mu(\lambda)^{-4}\chi_Y(\lambda)\phi(\v)$.
The group $G$ acts on $V_s$ by $(\tau_s(g)\phi)(\v)=\mu(\det g)\phi(^tg\v)$.
Let $\v_0\neq 0$ be a vector of $V$ such that the line
$\{\lambda\v_0;\lambda$ in $F\}$ is fixed under the action of $^tP$.
Explicitly, we take $\v_0=\,{}^t(0,0,0,1)$. It is clear that the map
$V_s\to W_s$, $\phi\mapsto\psi=\psi_\phi$, where
$\psi(g)=(\tau_s(g)\phi)(\v_0)=\mu(\det g)\phi(^tg\v_0)$,
is a $G$-module isomorphism (check that $\psi_{\tau_s(g)\phi}
=\pi_s(g)\psi_\phi$), with inverse $\psi\mapsto\phi=\phi_\psi$,
$\phi(\v)=\mu(\det g)^{-1}\psi(g)$ if $\v=\,{}^tg\v_0$ ($G$
acts transitively on $V-\{0\})$.

For $\v=\,{}^t(x,y,z,t)$ in $V$ put $\Vert\v\Vert=\max(\vert x\vert$,
$\vert y\vert$, $\vert z\vert$, $\vert t\vert)$. Let $V^0$ be the
quotient of the set $V^1$ of $\v$ in $V$ with $\Vert\v\Vert=1$, by the
equivalence relation $\v\sim\alpha\v$ if $\alpha$ is a
unit in $R$. Denote by $\Bbb P V$ the projective space
of lines in $V-\{0\}$. If $\Phi$ is a function on $V-\{0\}$
with $\Phi(\lambda\v)=\vert\lambda\vert^{-4}\Phi(\v)$ and
$d\v=dx~dy~dz~dt$, then $\Phi(\v)d\v$ is homogeneous of degree
zero. Define
$$\int_{\Bbb P V}\Phi(\v)d\v\quad\text{to be}\quad
\int_{V^0}\Phi(\v)d\v.$$
Clearly we have
$$\int_{\Bbb P V}\Phi(\v)d\v=\int_{\Bbb P V}
\Phi(g\v)d(g\v)=\vert\det g\vert\int_{\Bbb P V}
\Phi(g\v)d\v.$$
Put $\nu(x)=\vert x\vert$ and $m=2(s-1)$. Note that
$\nu/\mu_s=\mu_{-s}$. Put $\langle\w,\v\rangle=\,{}^t\w J\v$. Then
$\langle g\w,\theta(g)\v\rangle = \langle\w,\v\rangle$.

\proclaim {1. Proposition} The operator
$$T_s^Y:V_s\to V_{-s},\qquad (T_s^Y\phi)(\v)=\int_{\Bbb P V}
\phi(\w)\vert\langle\w,\v\rangle\vert^m \chi_Y(\langle\w,\v\rangle) d\w,$$
converges on $\Re s > 1/2$, and satisfies there
$T_s^Y\tau_s(g)=\tau_{-s}(\theta(g))T_s^Y$ for all
$g$ in $G$.\endproclaim

\demo {Proof} We have
$$\eqalign{(T_s^Y&(\tau_s(g)\phi))(\v)
=\smallint(\tau_s(g)\phi)(\w)\vert^t\w J\v\vert^m
\chi_Y(\langle\w,\v\rangle)d\w\cr
&=\mu(\det g)\smallint\phi(^tg\w)\vert^t\w J\v\vert^m
\chi_Y(\langle\w,\v\rangle)d\w\cr
&=\vert\det g\vert^{-1}\mu(\det g)\smallint\phi(\w)
\vert^t(^tg^{-1}\w)J\v\vert^m\chi_Y(\langle{}^tg^{-1}\w,\v\rangle)d\w\cr
&=(\mu/\nu)(\det g)\smallint\phi(\w)\vert^t\w J
\cdot J^{-1}g^{-1}J\v\vert^m\chi_Y(\langle\w,\theta(^tg)\v\rangle)d\w\cr
&=(\mu/\nu)(\det g)\smallint\phi(\w)\vert\langle\w,
\theta(^tg)\v\rangle\vert^m\chi_Y(\langle\w,\theta(^tg)\v\rangle)d\w\cr
&=(\nu/\mu)(\det\theta(g))\cdot
(T_s^Y\phi)({}^t\theta(g)\v)=[(\tau_{-s}(\theta (g)))(T_s^Y\phi)](\v)\cr}$$
for the functional equation.

For the convergence, we may assume that $\phi=1$ and ${}^t\v=(0,0,0,1)$,
so that the integral is $\int_R\,|x|^mdx$, which converges for $\Re m>-1$.
Our $m$ is $2s-2$, as required.
\enddemos

The spaces $V_s$ are isomorphic to the space $W$ of locally-constant
complex-valued functions $\phi$ on $V^1$ with $\phi(\lambda\v)
=\chi_Y(\lambda)\phi(\v)$ for all $\lambda\in R^\times$,
and $T_s^Y$ is equivalent to an operator $T^{Y,0}_s$ on $W$.
The proof of Proposition 1 implies also

\proclaim {1. Corollary} The operator $T^{Y,0}_s\circ\tau_s(g^{-1})$
is an integral operator with kernel
$$(\mu/\nu)(\det\theta (g))\vert\langle\w,\theta(^tg^{-1})\v\rangle
\vert^m\chi_Y(\langle\w,\theta(^tg^{-1})\v\rangle)
\qquad (\v,\w~\text{in}~V^1)$$
and trace
$$\tr[T^{Y,0}_s\circ\tau_s(g^{-1})]=(\nu/\mu)(\det g)
\int_{V^0}\vert^t\v gJ\v\vert^m\chi_Y(^t\v gJ\v) d\v.$$
\endproclaim

Next we normalize the operator $T^Y=T_s^{Y,0}$. Recalling that $\chi_Y$ is
unramified
($=1$ on $R^\times$, $\chi_Y(\pii)=-1$), we normalize $T^Y$ so that it acts
trivially on the one-dimensional space of $K$-fixed vectors in $V_s$.
This space is spanned by the function $\phi_0$ in $V_s$ with $\phi_0(\v)=1$
for all $\v$ in $V^0$. This is the only case studied in full in this paper.

Denote again by $\pii$ a generator of the maximal ideal of the ring
$R$ of integers in our local nonarchimedean field $F$ of odd residual
characteristic. Denote by $q$ the number of elements of the residue
field $R/\pii R$ of $R$. Normalize the absolute value by $|\pii|=q^{-1}$,
and the measures by $\Vol\{|x|\le 1\}=1$. Then $\Vol\{|x|=1\}=1-q^{-1}$,
and the volume of $V^0$ is $(1-q^{-4})/(1-q^{-1})=1+q^{-1}+q^{-2}+q^{-3}$.

\proclaim{2. Proposition} If $\v_0={}^t(0,0,0,1)$, we have
$$(T^Y\phi_0)(\v_0)={{1+q^{-2(s+1)}}\over{1+q^{1-2s}}}\phi_0(\v_0).$$
When $s=0$, the constant is $(1+q^{-2})(1+q)^{-1}$.
\endproclaim

\demo{Proof} Since $\chi_Y$ is unramified, we have
$$(T^Y\phi_0)(\v_0)=\int_{V^0}\phi_0(\v)|{}^t\v J\v_0|^m
\chi_Y(^t\v J\v_0)d\v=\int_{V^0}|x|^m \chi_Y(x) dxdydzdt$$
$$=\left[\int_{||\v||\le 1}-\int_{||\v||<1}\right]
|x|^m\chi_Y(x)dxdydzdt/\int_{|x|=1}dx$$
$$=(1+q^{-m-4})\int_{|x|\le 1}|x|^m\chi_Y(x) dx
/\int_{|x|=1}dx=(1+q^{-2(s+1)})/(1+q^{1-2s}),$$
since $m=2(s-1)$ and
$$\int_{|x|\le 1}|x|^m\chi_Y(x)dx=(1+q^{-m-1})^{-1}\int_{|x|=1}dx.$$
The proposition follows.
\enddemos

\heading Character computation for type I\endheading

For the $\theta$-conjugacy class of type I, represented by
$g=t\cdot\diag(\r,\s,\s,\r)$, the product
$${}^t\v gJ\v=(t,z,x,y)
\left(\matrix a_1\r & 0 & 0 & a_2D\r\\ 0 & b_1\s & b_2D\s & 0\\
0 & b_2\s & b_1\s& 0\\ a_2\r & 0 & 0 & a_1\r\endmatrix\right)
\left(\matrix 0 & 0 & 0 & 1\\ 0 & 0 & 1 & 0\\
0 &-1 & 0 & 0\\-1 & 0 & 0 & 0\endmatrix\right)
\left(\matrix t\\ z\\ x\\ y\endmatrix\right)$$
is equal to
$$-t^2a_2D\r-z^2b_2D\s+x^2b_2\s+y^2a_2\r.$$
Note that $\r$ and $\s$ range over a set of representatives for
$F^\times/N_{E/F}E^{\times}$.

By Corollary 1, we need to compute $$({\nu\over{\mu}})(\det g)
{\Delta(g\theta)\over{\Delta_C(Ng)}} \int_{V^0}|{}^t\v gJ\v|^m
\chi_Y({}^t\v gJ\v)d\v$$
$$={|\r\s|^{1-s}|4a_2b_2D|\over{|(a_1^2-a_2^2D)(b_1^2-b_2^2D)|^{s/2}}}
\int_{V^0}|\alpha|^{2(s-1)}\chi_Y(\alpha) dxdydzdt.$$
Here $\alpha$ is $x^2b_2\s+y^2a_2\r-z^2b_2D\s-t^2a_2D\r$.
Put $\r'=-{a_2\over{b_2}}{\r\over\s}$. Thus we need to compute
the value at $s=0$ of the product of
$$\chi_Y(b_2\s)|{\r\over{\s}}|^{-s}|4D\r'||(({a_1\over{b_2}})^2
-({a_2\over{b_2}})^2D)(({b_1\over{b_2}})^2-D)|^{-s/2}$$
with the integral $I_s^Y(\r',D)$, where $Q=x^2-\r y^2-D z^2+\r D t^2$ and
$$I_s^Y(\r,D)=\int_{V^0}|Q|^{2(s-1)}\chi_Y(Q)dxdydzdt.$$

\proclaim{I. Theorem}
When $Y/F$ is unramified, the value of $I_s^Y(\r,D)$ at $s=0$ is $0$.
\endproclaim

\demo{Proof} Consider the case when the quadratic form
$Q=x^2-\r y^2-D z^2+\r D t^2$ is anisotropic (does not represent zero).
Thus $D=\pii$ and $\r\in R^\times-R^{\times 2}$ (hence $|\r|=1$,
$|D|=1/q$), or $D\in R^\times-R^{\times 2}$ and $\r=\pii$. The second
case being equivalent to the first, it suffices to deal with the first
case.

The domain $\max\{|x|,|y|,|z|,|t|\}=1$ is the disjoint union of
$\{|x|=1\}$, $\{|x|<1,|y|=1\}$, $\{|x|<1,|y|<1,|z|=1\}$
and $\{|x|<1,|y|<1,|z|<1,|t|=1\}$.
Note that $\chi_Y(x^2-\r y^2-D z^2+\r D t^2)$ is equal to $1$ on
the first two subdomains and equals $-1$ on the other two.
Thus the integral $I_s^Y(\r,D)$ is the quotient by $\int_{|x|=1}dx$ of
$$\int_{|x|=1} dx+\int\int_{|x|<1,|y|=1} dx dy
-q^{-m}\int\int\int_{|x|<1,|y|<1,|z|=1} dx dy dz$$
$$-q^{-m}\int\int\int\int_{|x|<1,|y|<1,|z|<1,|t|=1}dxdydzdt$$
$$=1+q^{-1}-q^{-m-2}-q^{-m-3}
=1+q^{-1}-q^{-2s}-q^{-2s-1}.$$
The value at $s=0$ is $0$ and thus
the theorem follows when the quadratic form is
anisotropic.
\enddemo

We then turn to the case when the quadratic form is isotropic.
Recall that $\r$ ranges over a set of representatives for
$F^\times/N_{E/F}E^\times$, $E=F(\sqrt D)$.
Thus $D\in F-F^2$, and we may assume that $|D|$ and $|\r|$
lie in $\{1,q^{-1}\}$.

\proclaim{I.1. Proposition \fz4}
When the quadratic form
$x^2-\r y^2-D z^2+\r D t^2$ is isotropic, $\r$ lies in $N_{E/F}E^\times$,
and we may assume that the quadratic form takes one of three shapes:
$$x^2-y^2-D z^2+D t^2,\,\,\, D\in R^\times-R^{\times 2};\quad
x^2+\pii y^2-\pii z^2-\pii^2t^2;\quad x^2-y^2-\pii z^2+\pii t^2.$$
\endproclaim

The set $V^0=V/\sim$, where
$V=\{\v=(x,y,z,t)\in R^4;\max\{|x|,|y|,|z|,|t|\}=1\}$
and $\sim$ is the equivalence relation $\v\sim\alpha\v$ for
$\alpha\in R^\times$, is the disjoint union of the subsets
$$V_n^0=V_n^0(\r,D)=V_n(\r,D)/\sim,$$ where
$$V_n=V_n(\r,D)=\{\v;\max\{|x|,|y|,|z|,|t|\}=1,
|x^2-\r y^2-D z^2+\r D t^2|=1/q^n\},$$
over $n\ge 0$, and of $\{\v; x^2-\r y^2-D z^2+\r D t^2=0\}/\sim$,
a set of measure zero.

Thus the integral $I_s^Y(\r,D)$ coincides with the sum
$$\sum_{n=0}^{\infty}(-1)^nq^{-nm}\Vol(V_n^0(\r,D)).$$
When the quadratic form represents zero the problem
is then to compute the volumes
$$\Vol(V_n^0(\r,D))=\Vol(V_n(\r,D))/(1-1/q)\qquad (n\ge 0).$$

We need some results from [FZ4]:

\proclaim{I.0. Lemma \fz4}
When $c^2\in R^{\times 2}$ and $n\ge 1$, we have
$$\int_{|c^2-x^2|=q^{-n}}dx~=~\frac{2}{q^n}\left(1-\frac{1}{q}\right).$$
\endproclaim

\proclaim{I.1. Lemma \fz4}
When $D=\pii$ and $\r=1$, thus $|\r D|=1/q$, we have
$$\Vol(V_n^0)=\cases
1-1/q,               &\qquad\text{if}\ n=0,\\
q^{-1}(1-1/q)(2+1/q) , &\qquad\text{if}\ n=1,\\
2q^{-n}(1-1/q)(1+1/q), &\qquad\text{if}\ n\ge 2.
\endcases$$
\endproclaim

\proclaim{I.2. Lemma \fz4}
When $D=\pii$ and $\r=-\pii$, thus $|\r D|=1/q^2$, we have
$$\Vol(V_n^0)=\cases
1,                    &\qquad\text{if}\ n=0,\\
q^{-1}(1-1/q),        &\qquad\text{if}\ n=1,\\
q^{-2}(2-1/q-2/q^2),  &\qquad\text{if}\ n=2,\\
2q^{-n}(1-1/q)(1+1/q),&\qquad\text{if}\ n\ge 3.
\endcases $$
\endproclaim

\proclaim{I.3. Lemma \fz4}
When $E/F$ is unramified, thus $|\r D|=1$, we have
$$\Vol(V_n^0)=\cases
1-1/q^2,                    &\qquad\text{if}\ n=0,\\
q^{-n}(1-1/q)(1+2/q+1/q^2), &\qquad\text{if}\ n\ge 1.
\endcases$$
\endproclaim

\demo{Proof of Theorem I} We are now ready to complete the
proof of Theorem I in the isotropic case. Recall that
we need to compute the value at $s=0$ ($m=-2$) of
$I_s^Y(\r,D)$. Here $I_s^Y(\r,D)$ coincides with the sum
$$\sum_{n=0}^{\infty}(-1)^nq^{-nm}\Vol(V_n^0(\r,D))$$
which converges for $m>-1$ by Proposition 1 or alternatively
by Lemmas I.1-I.3. The value at $m=-2$ is obtained then by
analytic continuation of this sum.

\n {\it Case of Lemma} I.1. The integral $I_s^Y(\r,D)$ is equal to
$$\Vol(V_0^0)-q^{-m}\Vol(V_1^0)+\sum_{n=2}^{\infty}(-1)^nq^{-nm}\Vol(V_n^0)$$
$$=1-\frac{1}{q}-\frac{1}{q}\left(1-\frac{1}{q}\right)
\left(2+\frac{1}{q}\right)
\frac{1}{q^m}+2\left(1-\frac{1}{q}\right)\left(1+\frac{1}{q}\right)
q^{-2(m+1)}\left(1+\frac{1}{q^{m+1}}\right)^{-1}.$$
When $m=-2$, this is
$$1-\frac{1}{q}-q\left(2-\frac{1}{q}-\frac{1}{q^2}\right)
+2\left(1-\frac{1}{q}\right)\left(1+\frac{1}{q}\right)\frac{q^2}{1+q}
=0.$$

\n {\it Case of Lemma} I.2. The integral $I_s^Y(\r,D)$ is equal to
$$\Vol(V_0^0)-q^{-m}\Vol(V_1^0)+q^{-2m}\Vol(V_2^0)
+\sum_{n=3}^{\infty}(-1)^nq^{-nm}\Vol(V_n^0)$$
$$=1-\frac{1}{q}\left(1-\frac{1}{q}\right)q^{-m}
+\frac{1}{q^2}\left(2-\frac{1}{q}-\frac{2}{q^2}\right)q^{-2m}$$
$$-2\left(1-\frac{1}{q}\right)\left(1+\frac{1}{q}\right)
q^{-3(m+1)}\left(1+\frac{1}{q^{m+1}}\right)^{-1}.$$
When $m=-2$, this is
$$1-\frac{1}{q}\left(1-\frac{1}{q}\right)q^2
+\frac{1}{q^2}\left(2-\frac{1}{q}-\frac{2}{q^2}\right)q^4
- 2\left(1-\frac{1}{q}\right)\left(1+\frac{1}{q}\right)
\frac{q^3}{1+q}.$$
Once simplified this is equal to $0$.

\n {\it Case of Lemma} I.3. The integral $I_s^Y(\r,D)$ is equal to
$$\Vol(V_0^0)+\sum_{n=1}^{\infty}(-1)^nq^{-nm}\Vol(V_n^0)$$
$$=1-\frac{1}{q^2}-\left(1-\frac{1}{q}\right)
\left(1+\frac{2}{q}+\frac{1}{q^2}\right)
q^{-(m+1)}\left(1+\frac{1}{q^{m+1}}\right)^{-1}.$$
When $m=-2$, this is
$$=1-\frac{1}{q^2}-\left(1-\frac{1}{q}\right)
\left(1+\frac{1}{q}\right)^2\frac{q}{1+q}=0.$$
The theorem follows.
\enddemos

\heading Character computation for type II\endheading

For the $\theta$-conjugacy class of type II, represented by
$g=t\cdot\diag(\r,\s,\s,\r)$, the product
$${}^t\v gJ\v=(t,z,x,y)
\left(\matrix a_1\r & 0 & 0 & a_2D\r\\ 0 & b_1\s & b_2AD\s & 0\\
0 & b_2\s & b_1\s& 0\\ a_2\r & 0 & 0 & a_1\r\endmatrix\right)
\left(\matrix 0 & 0 & 0 & 1\\ 0 & 0 & 1 & 0\\
0 &-1 & 0 & 0\\-1 & 0 & 0 & 0\endmatrix\right)
\left(\matrix t\\ z\\ x\\ y\endmatrix\right)$$
is equal to
$$-t^2a_2D\r-z^2b_2AD\s+x^2b_2\s+y^2a_2\r=b_2\s(x^2-y^2\r'-z^2AD+t^2D\r').$$
Here $a_1+a_2\sqrt D\in E_1^\times$ ($E_1=F(\sqrt{D})$) and
$b_1+b_2\sqrt{AD}\in E_2^\times$ ($E_2=F(\sqrt{AD})$), and $\r'=-a_2\r/b_2\s$.
As $\r$ ranges over a set of representatives for
$F^\times/N_{E_1/F}E_1^{\times}$ (and $\s$ for
$F^\times/N_{E_2/F}E_2^{\times}$), we may rename $\r'$ by $\r$.

Thus, by Corollary 1, we need to compute the value at $s=0$ of the product of
$$\chi_Y(b_2\s)|{\r\over{\s}}|^{-s}|4\r'D\sqrt{A}|
|(({a_1\over{b_2}})^2-({a_2\over{b_2}})^2D)(({b_1\over{b_2}})^2-AD)|^{-s/2}$$
and the value when $\r$ is $\r'$ and $Q=x^2-\r y^2-AD z^2+\r D t^2$ of the
integral
$$I_s^Y(\r,A,D)=\int_{V^0}|Q|^{2(s-1)}
\chi_Y(Q)dxdydzdt.$$

The property of the numbers $A$, $D$ and $AD$ that we need is that
their square roots generate the three distinct quadratic extensions
of $F$. Thus we may assume that $\{A,D,AD\}=\{u,\pii,u\pii\}$, where
$u\in R^\times-R^{\times 2}$. Of course with this normalization $AD$
is no longer the product of $A$ and $D$, but its representative in the
set $\{1,u,\pii,u\pii\}$ mod $F^{\times 2}$. Since $\r$ ranges over a
set of representatives for $F^\times/N_{E_1/F}E_1^\times$, it can be
assumed to range over $\{1,\pii\}$ if $D=u$, and over $\{1,u\}$
if $|D|=|\pii|$.

In this section we prove

\proclaim{II. Theorem} When $Y/F$ is unramified, the value of
$$\chi_Y(b_2\s)|4\r D\sqrt{A}|I_s^Y(\r,A,D)/(T^Y\phi_0)(\v_0)$$
at $s=0$ is $-2\chi_Y(b_2\s)\delta(Y,E_3)$.
\endproclaim

Recall that $E_3=F(\sqrt{A})$. As usual, $\delta(Y,E_3)$ is 1 if $Y=E_3$
and 0 if $Y\not=E_3$.

The meaning of this result is that the twisted character of $\pi_Y$ on
elements of tori of type II relates to values of the trivial character
on $\uC_Y(F)$, $Y=E_3$, on the torus which splits over $E$. It does not
relate to such values on $\uC_{Y'}(F)$, $Y'\not=E_3$.

Recall Lemma II.0 from [FZ4].

\proclaim{II.0. Lemma \fz4} The quadratic form
$x^2-\r y^2-AD z^2+\r D t^2$ takes one of six forms$:$\hb
$x^2-y^2+\pii(t^2-uz^2)$,
$x^2-u y^2+u\pii(t^2-z^2)$,
$x^2-y^2+ut^2-u\pii z^2$,
$x^2-y^2-u z^2+\pii t^2$,
$x^2-\pii y^2+u\pii(t^2-z^2)$,
$x^2-u y^2-u z^2+u\pii t^2$,
where $u\in R^\times-R^{\times 2}$.
It is always isotropic.
\endproclaim

The set $V^0=V/\sim$, where
$V=\{\v=(x,y,z,t)\in R^4;\max\{|x|,|y|,|z|,|t|\}=1\}$
and $\sim$ is the equivalence relation $\v\sim\alpha\v$ for
$\alpha\in R^\times$, is the disjoint union of the subsets
$$V_n^0=V_n^0(\r,A,D)=V_n(\r,A,D)/\sim,$$ where
$$V_n=V_n(\r,A,D)=\{\v;\max\{|x|,|y|,|z|,|t|\}=1,
|x^2-\r y^2-AD z^2+\r D t^2|=1/q^n\},$$
over $n\ge 0$, and of $\{\v; x^2-\r y^2-AD z^2+\r D t^2=0\}/\sim$,
a set of measure zero.

Since $Y/F$ is unramified, the integral $I_s^Y(\r,A,D)$ is
equal to
$$\sum_{n=0}^{\infty}(-1)^nq^{-nm}\Vol(V_n^0(\r,A,D)).$$
The problem is then to compute the volumes
$$\Vol(V_n^0(\r,A,D))=\Vol(V_n(\r,A,D))/(1-1/q)\qquad (n\ge 0).$$

Choose $u$ to be a non square unit. To prove Theorem II, by Lemma II.0
we need precisely the following Lemmas from [FZ4]. In Lemmas II.1 and II.2, 
$E_3=F(\sqrt{A})$ is unramified over $F$.

\proclaim{II.1. Lemma \fz4} When the quadratic form is $x^2-y^2+\pii (t^2-uz^2)$
$($thus $\r=1$, $A=u$, $D=\pii$ up to squares$)$, we have
$$\Vol(V_n^0)=\cases
1-1/q,           &\qquad\text{if}\ n=0,\\
2/q-1/q^2+1/q^3, &\qquad\text{if}\ n=1,\\
2q^{-n}(1-1/q),  &\qquad\text{if}\ n\ge 2.
\endcases$$
\endproclaim

\proclaim{II.2. Lemma \fz4}
When the quadratic form is $x^2-uy^2+u\pii (t^2-z^2)$
$($thus $\r=u$, $A=u$, $D=\pii$ up to squares$)$, we have
$$\Vol(V_n^0)=\cases
1+1/q,             &\qquad\text{if}\ n=0,\\
q^{-2}(1-1/q),     &\qquad\text{if}\ n=1,\\
2q^{-(n+1)}(1-1/q),&\qquad\text{if}\ n\ge 2.
\endcases$$
\endproclaim

\proclaim{II.3. Lemma \fz4}
When the quadratic form is $x^2-y^2+ut^2-u\pii z^2$
or $x^2-y^2-uz^2+\pii t^2$, we have
$$\Vol(V_n^0)=\cases
1,              &\qquad\text{if}\ n=0,\\
1/q,            &\qquad\text{if}\ n=1,\\
q^{-n}(1-1/q^2),&\qquad\text{if}\ n\ge 2.
\endcases$$
\endproclaim

\proclaim{II.4. Lemma \fz4}
When the quadratic form is
$x^2-\pii y^2+u\pii (t^2-z^2)$, we have
$$\Vol(V_n^0)=\cases
1,              &\qquad\text{if}\ n=0,\\
1/q,            &\qquad\text{if}\ n=1,\\
q^{-n}(1-1/q^2),&\qquad\text{if}\ n\ge 2.
\endcases$$
\endproclaim

\proclaim{II.5. Lemma \fz4}
When the quadratic form is $x^2-uy^2-uz^2+u\pii t^2$, we have
$$\Vol(V_n^0)=\cases
1,              &\qquad\text{if}\ n=0,\\
1/q,            &\qquad\text{if}\ n=1,\\
q^{-n}(1-1/q^2),&\qquad\text{if}\ n\ge 2.
\endcases$$
\endproclaim

\demo{Proof of Theorem II} To prove Theorem II, recall that we need 
to compute the value at $s=0$ ($m=-2$) of the product
$$\chi_Y(b_2\s)|4\r D\sqrt{A}|I_s^Y(\r,A,D)/(T^Y\phi_0)(\v_0).$$
Here $I_s^Y(\r,A,D)$ is equal to the sum
$$\sum_{n=0}^{\infty}(-1)^nq^{-nm}\Vol(V_n^0(\r,D))$$
which converges for $m>-1$ by Proposition 1 or alternatively
by Lemmas II.1-II.3. The value at $m=-2$ is obtained then by
analytic continuation of this sum.

\n {\it Case of Lemma} II.1. We have $|4\r D\sqrt{A}|=1/q$, and the integral
$I_s^Y(\r,A,D)$ is equal to
$$\Vol(V_0^0)-q^{-m}\Vol(V_1^0)+\sum_{n=2}^{\infty}(-1)^nq^{-nm}\Vol(V_n^0)$$
$$=1-\frac{1}{q}-\left(\frac{2}{q}-\frac{1}{q^2}+\frac{1}{q^3}\right)
\frac{1}{q^m}+2\left(1-\frac{1}{q}\right)q^{-2(m+1)}
\left(1+\frac{1}{q^{m+1}}\right)^{-1}.$$
When $m=-2$, this is
$$1-\frac{1}{q}-q^2\left(\frac{2}{q}-\frac{1}{q^2}+\frac{1}{q^3}\right)
+2\left(1-\frac{1}{q}\right)\frac{q^2}{1+q}
=\frac{-2q}{1+q}\left(1+\frac{1}{q^2}\right).$$
Multiplying by $|4\r D\sqrt{A}|=1/q$ we obtain $-2(1+1/q^2)(1+q)^{-1}$.
We are done by Proposition 2.

\n {\it Case of Lemma} II.2. We have $|4\r D\sqrt{A}|=1/q$, and the integral
$I_s^Y(\r,A,D)$ is equal to
$$\Vol(V_0^0)-q^{-m}\Vol(V_1^0)+\sum_{n=2}^{\infty}(-1)^nq^{-nm}\Vol(V_n^0)$$
$$=1+\frac{1}{q}-\frac{1}{q^2}\left(1-\frac{1}{q}\right)
\frac{1}{q^m}+\frac{2}{q}\left(1-\frac{1}{q}\right)q^{-2(m+1)}
\left(1+\frac{1}{q^{m+1}}\right)^{-1}.$$
When $m=-2$, this is
$$1+\frac{1}{q}-1+\frac{1}{q}+\frac{2}{q}\left(1-\frac{1}{q}\right)
\frac{q^2}{1+q}= \frac{2q}{1+q}\left(1+\frac{1}{q^2}\right).$$
Multiplying by $|4\r D\sqrt{A}|=1/q$ we obtain $2(1+1/q^2)(1+q)^{-1}$.

\n {\it Case of Lemmas} II.3, II.4, II.5. The integral $I_s^Y(\r,A,D)$
is equal to
$$\Vol(V_0^0)-q^{-m}\Vol(V_1^0)+\sum_{n=2}^{\infty}(-1)^nq^{-nm}\Vol(V_n^0)$$
$$=1-\frac{1}{q}\frac{1}{q^m}+\left(1-\frac{1}{q}\right)
\left(1+\frac{1}{q}\right)q^{-2(m+1)}\left(1+\frac{1}{q^{m+1}}\right)^{-1}.$$
When $m=-2$, this is
$$1-\frac{1}{q}q^2+\left(1-\frac{1}{q}\right)
\left(1+\frac{1}{q}\right)\frac{q^2}{1+q}=0.$$
\enddemos

\heading Character computation for type III\endheading

For the $\theta$-conjugacy class of type III we write out the
representative $g=t\cdot\diag({\pmb r},{\pmb r})$ as
$$
\left(\matrix
a_1r_1+a_2r_2A & (a_1r_2+a_2r_1)A & (b_1r_1+b_2r_2A)D & (b_1r_2+b_2r_1)AD\\
a_1r_2+a_2r_1  & a_1r_1+a_2r_2A   & (b_1r_2+b_2r_1)D  & (b_1r_1+b_2r_2A)D\\
b_1r_1+b_2r_2A & (b_1r_2+b_2r_1)A &  a_1r_1+a_2r_2A   & (a_1r_2+a_2r_1)A\\
b_1r_2+b_2r_1  &  b_1r_1+b_2r_2A  &  a_1r_2+a_2r_1    &  a_1r_1+a_2r_2A
\endmatrix\right).$$
The product ${}^t\v gJ\v$ (where ${}^t\v=(x,y,z,t)$) is equal to
$$(b_1r_2+b_2r_1)(t^2+z^2A-y^2D-x^2AD)
+2(b_1r_1+b_2r_2A)(zt-xyD),$$
where $a_1+a_2\sqrt A\in E_3^\times$ and $b_1+b_2\sqrt A\in E_3^\times$.
The trace is a function of $g$, and
$r=r_1+r_2\sqrt{A}$ ranges over a set of representatives in $E_3^\times$
($E_3=F(\sqrt{A})$) for $E_3^\times/N_{E/E_3}E^{\times}$.

Define the quadratic form $Q=Q(x,y,z,t)$ to be
$${{rb-\tau(rb)}\over{2\sqrt A}}(t^2+z^2A-y^2D-x^2AD)+(rb+\tau(rb))(zt-xyD).$$
Set $I_s^Y(r,A,D)$ to be equal to
$$\int_{V^0}|Q|^{2(s-1)}\chi_Y(Q) dxdydzdt.$$
The property of the numbers $A$, $D$ and $AD$ that we need is that
their square roots generate the three distinct quadratic extensions
of $F$. Thus we may assume that $\{A,D,AD\}=\{u,\pii,u\pii\}$, where
$u\in R^\times-R^{\times 2}$. Of course with this normalization $AD$
is no longer the product of $A$ and $D$, but its representative in the
set $\{1,u,\pii,u\pii\}$ mod $F^{\times 2}$.

\proclaim{III.1. Proposition} $(i)$ If $D=u$ and $A=\pii$ $($or $\pii u)$
then $\sqrt{A}\not\in N_{E/E_3}E^\times=A^\Z R_3^\times$.\hb
$(ii)$ If $A=u$ and $-1\in R^{\times 2}$, and $D=\pii$ $($or $\pii u)$ then
$\sqrt{A}\not\in N_{E/E_3}E^\times=(-D)^\Z R_3^{\times 2}$.\hb
$(iii)$ If $A=u=-1\not\in R^{\times 2}$ and $D=\pii$ $($or $\pii u)$ then
there is $d\in R^\times$ with $d^2+1\in -R^{\times 2}=R^\times-R^{\times 2}$,
hence $d+i\in R_3^{\times}-R_3^{\times 2}$ $(i=\sqrt{A})$ and so $d+i\in
E_3^\times-N_{E/E_3}E^\times$.
\endproclaim

\demo{Proof} For (iii) note that $R^\times/\{1+\pii R\}$ is the
multiplicative group of a finite field $\F$ of $q$ elements. There
are $1+{1\over 2}(q-1)$ elements in each of the sets $\{1+x^2;\,x\in\F\}$
and $\{-y^2;\,y\in\F\}$. As $2(1+{1\over 2}(q-1))>q$, there are $x$, $y$
with $1+x^2=-y^2$. But $y\not=0$ as $-1\not\in\F^{\times 2}$. Hence there
is $x$ with $1+x^2\not\in\F^{\times 2}$, and our $d$ exists.
\enddemos

Since $r$ ranges over a set of representatives for
$E_3^\times/N_{E/E_3}E^{\times}$, by Proposition III.1
we can choose $br$ to be $1$ or $\sqrt A$ or $d+i$.
Correspondingly the quadratic form takes one of the three shapes
$$t^2+z^2A-y^2D-x^2AD,\qquad zt-xyD,\qquad\text{or}\qquad
t^2-z^2-y^2D+x^2D+2d(zt-xyD).$$

\proclaim{III. Theorem} When $Y/F$ is unramified,
the value of $I_s^Y(\r,A,D)$ at $s=0$ is $0$.
\endproclaim

\demo{Proof} Assume that $br=\sqrt A\not\in N_{E/E_3}E^\times$, thus
$|br\tau(br)D|=|AD|$, and the quadratic form is $t^2+z^2A-y^2D-x^2AD$.
If $|A|=1/q$ or $-1$ is a square, we can replace $A$ with $-A$. The
quadratic form then becomes the same as that of type I. The result of
the computation is $0$, see proof of Theorem I, case of anisotropic
quadratic forms and we are done in this case.

If $A=-1$, $br=d+i\not\in N_{E/E_3}E^\times$, the quadratic form is
$t^2-z^2-y^2D+x^2D+2d(zt-xyD)$. It is equal to $X^2-uY^2-D(Z^2-uT^2)$
with $X=t+dz$, $Y=z$, $Z=y+dx$, $T=x$ and $u=d^2+1
\in R^\times-R^{\times 2}$. Since $|D|=1/q$ the quadratic form is
anisotropic and the result of the computation is $0$ by the proof of
Theorem I, case of anisotropic quadratic forms.

Assume that $br=1$, thus $|br\tau(br)D|=|D|$ and the
quadratic form is $zt-xyD$. Then it is ${1\over 4}$ times
$(z+t)^2-(z-t)^2-D[(x+y)^2-(x-y)^2].$
Since $\max\{|x|,|y|,|z|,|t|\}=1$ implies
$\max\{|x+y|,|x-y|,|z+t|,|z-t|\}=1$, the result of the computation
is $0$ by the proof of Theorem I, cases of Lemmas I.1 and I.3.
The theorem follows.
\enddemos

\heading Character computation for type IV\endheading

For the $\theta$-conjugacy class of type IV we write the representative
$g=t\cdot\diag({\pmb r},{\pmb r})$ (where $t=h^{-1}t^\ast h$, $t^\ast=
\diag(\alpha,\sigma\alpha,\sigma^3\alpha,\sigma^2\alpha)$) as
$$
\left(\matrix
a_1r_1+a_2r_2A & (a_1r_2+a_2r_1)A & (b'_1r_1+b'_2r_2A)D &(b'_1r_2+b'_2r_1)AD\\
a_1r_2+a_2r_1  & a_1r_1+a_2r_2A   & (b'_1r_2+b'_2r_1)D  &(b'_1r_1+b'_2r_2A)D\\
b_1r_1+b_2r_2A & (b_1r_2+b_2r_1)A &  a_1r_1+a_2r_2A     & (a_1r_2+a_2r_1)A\\
b_1r_2+b_2r_1  &  b_1r_1+b_2r_2A  &  a_1r_2+a_2r_1      &  a_1r_1+a_2r_2A
\endmatrix\right).$$
Here $E_3=F(\sqrt{A})$ is a quadratic extension of $F$ and $E=E_3(\sqrt D)$
is a quadratic extension of $E_3$, thus $A\in F-F^2$ and $D=d_1+d_2\sqrt A\in
E_3-E_3^2$, $d_i\in F$.

If $-1\in F^{\times 2}$ we can and do take $D=\sqrt{A}$, where $A$ is a
nonsquare unit $u$ if $E_3/E$ is unramified, or a uniformizer $\pii$ if
$E_3/F$ is ramified. If $-1\not\in F^{\times 2}$ and $E_3/F$ is ramified,
once again we may and do take $A=\pii$ and $D=\sqrt{A}$.

If $-1\not\in F^{\times 2}$ and $E_3/F$ is unramified, take $A=-1$ and
note that a primitive 4th root $\zeta=i$ of 1 lies in $E_3$ (and generates
it over $F$). Then $E/E_3$ is unramified, generated by $\sqrt{D}$,
$D=d_1+id_2$, and we can (and do) take $d_2=1$ and a unit $d_1=d$ in
$F^\times$ such that $d^2+1\not\in F^{\times 2}$. Then $D=d+i\not\in
E_3^{\times 2}$. The existence of $d$ is shown as in the proof of
Proposition III.1.

Further $\alpha=a+b\sqrt{D}\in E^\times$, where
$a=a_1+a_2\sqrt A\in E_3^\times$, $b=b_1+b_2\sqrt A\in E_3^\times$, and
$r=r_1+r_2\sqrt{A}\in E_3^\times/N_{E/E_3}E^\times$. The relation
$bD=b'_1+b'_2\sqrt A$ defines $b'_1=b_1d_1+b_2d_2A$ and $b'_2=b_2d_1+b_1d_2$.

Recall from Corollary 1 that we need to compute
$$\left({\nu\over{\mu}}\right)(\det g){\Delta(g\theta)
\over{\Delta_C(Ng)}}\int_{V^0}|{}^t\v gJ\v|^m\chi_Y({}^t\v gJ\v)d\v.\eqno(\ast)$$
Since $\det g=\alpha r\cdot\sigma(\alpha r)\cdot
\sigma^3(\alpha r)\cdot\sigma^2(\alpha r)$, we have
$$\left({\nu\over{\mu}}\right)(\det g){\Delta(g\theta)\over{\Delta_C(Ng)}}
=|\det g|^{(1-s)/2}\left|{(\alpha r-\sigma^2(\alpha r))^2\over{\alpha
r\sigma^2(\alpha r)}}\cdot {\sigma(\alpha r-\sigma^2(\alpha r))^2\over
{\sigma(\alpha r)\sigma^3(\alpha r)}}\right|^{1/2}$$
$$={|4brD\sigma(brD)|\over |r^2(a^2-b^2D)
\sigma(r^2(a^2-b^2D))|^{s/2}}.$$
When $s=0$, this is $|brD \sigma(brD)|$.

The product ${}^t\v gJ\v$ (where ${}^t\v=(x,y,z,t)$) is then equal to
$$(b_1r_2+b_2r_1)(t^2+z^2A)-(b'_1r_2+b'_2r_1)(y^2+x^2A)
+2(b_1r_1+b_2r_2A)zt-2(b'_1r_1+b'_2r_2A)xy.$$

Since $bD=b'_1+b'_2\sqrt A$, this is
$${{br-\sigma(br)}\over{2\sqrt A}}(t^2+z^2A)+(br+\sigma(br))zt$$
$$-{{brD-\sigma(brD)}\over{2\sqrt A}}(y^2+x^2A)-(brD+\sigma(brD))xy.$$

Note that $r$ ranges over a set of representatives for
$E_3^\times/N_{E/E_3}E^\times$, and $b$ lies in $E_3^\times$.
As $b$ is fixed, we may take $br$ to range over
$E_3^\times/N_{E/E_3}E^\times$.

Further, note that $E_3/F$ is unramified if and only if $E/E_3$ is
unramified. Hence $br$ can be taken to range over $\{1,\pii\}$ if $E_3/F$
is unramified, and over $\{1,u\}$ if $E_3/F$ is ramified, where $\pii$ is a
uniformizer in $F$ and $u$ is a nonsquare unit in $F$, in these two cases.
Thus in both cases we have that $\sigma(br)=br$, and the quadratic form is
equal to $brQ$, where
$$Q=Q(x,y,z,t)=2zt-{{D-\sigma(D)}\over{2\sqrt A}}(y^2+x^2A)-(D+\sigma(D))xy.$$
Thus we need to compute the value at $s=0$ of the product of
$|brD\sigma(brD)|$, $\chi_Y(br)|br|^{2(s-1)}$ and
$$I_s^Y(\r,A,D)=\int_{V^0}|Q|^{2(s-1)}\chi_Y(Q)dxdydzdt.$$

\proclaim{IV. Theorem} When $Y/F$ is unramified, the the value of
$$\chi_Y(br)|br|^{2(s-1)}|brD\sigma(brD)|I_s^Y(r,A,D)/(T^Y\phi_0)(\v_0)$$
at $s=0$ is $-2\chi_Y(br)\delta(Y,E_3)$.
\endproclaim

To prove this theorem we need some results from [FZ4].

\proclaim{IV.1. Proposition \fz4}
Up to a change of coordinates, the quadratic form
$$2zt-{{D-\sigma(D)}\over{2\sqrt A}}(y^2+x^2A)-(D+\sigma(D))xy$$
is equal to either $x^2+\pii y^2-2zt$ or $x^2-uy^2-2zt$
with $u\in R^\times-R^{\times 2}$. It is always isotropic.
\endproclaim

Recall that $Y/F$ is unramified. Then the integral
$I_s^Y(\r,A,D)$ is equal to
$$\sum_{n=0}^{\infty}(-1)^nq^{-nm}\Vol(V_n^0(\r,A,D)).$$
The problem is then to compute the volumes
$$\Vol(V_n^0(\r,A,D))=\Vol(V_n(\r,A,D))/(1-1/q)\qquad (n\ge 0).$$

\proclaim{IV.2. Lemma \fz4}
When the quadratic form is $x^2-uy^2-2zt$, we have
$$\Vol(V_n^0)=\cases
1+1/q^2,                &\qquad\text{if}\ n=0,\\
q^{-n}(1-1/q)(1+1/q^2), &\qquad\text{if}\ n\ge 1.
\endcases$$
\endproclaim

\demo{Proof of Theorem IV} To prove Theorem IV, recall that we need 
to compute the value at $s=0$ ($m=-2$) of $I_s^Y(\r,A,D)$. Since $Y/F$ 
is unramified, the integral $I_s^Y(\r,A,D)$ coincides
with the sum $$\sum_{n=0}^{\infty}(-1)^nq^{-nm}\Vol(V_n^0(A,D))$$
which converges for $m>-1$. The value at $m=-2$ is obtained then by
analytic continuation of this sum.
\medskip
\n {\it Case of $x^2+\pii y^2-2zt$}. Make a change of variables
$z\mapsto 2u^{-1} z'$, followed by $z'\mapsto z$. Thus the quadratic
form is equal to
$$-u^{-1}((z-t)^2-(z+t)^2-u x^2-u\pii y^2).$$
Note that up to a multiple by a unit, this is a form of Lemma II.3.
Since $\max\{|z|,|t|\}=1$ implies $\max\{|z+t|,|z-t|\}=1$,
the result of that lemma holds for our quadratic form as well.
In this case $E_3/F$ is ramified, and our integral is zero.
\medskip
\n {\it Case of $x^2-uy^2-2zt$}. This is the case where $E_3/F$ is
unramified. By Lemma IV.2, the integral
$$I_s^Y(\r,A,D) =\Vol(V_0^0)+\sum_{n=1}^{\infty}(-1)^nq^{-nm}\Vol(V_n^0)$$
is equal to
$$1+\frac{1}{q^2}+\left(1-\frac{1}{q}\right)
\left(1+\frac{1}{q^2}\right)\frac{-1}{q^{(m+1)}}
\left(1-\frac{-1}{q^{m+1}}\right)^{-1}.$$
When $m=-2$, this is
$$1+\frac{1}{q^2}+\left(1-\frac{1}{q}\right)\left(1+\frac{1}{q^2}\right)
\frac{-q}{q+1}=\frac{2}{1+q}\left(1+\frac{1}{q^2}\right).$$
The theorem follows by Proposition 2.
\enddemos
\medskip

\def\refe#1#2{\n\hangindent 5em\hangafter1\hbox to 5em{\hfil#1\quad}#2}
\subheading{References}
\medskip
\refe{[BZ]}{I. Bernstein, A. Zelevinsky, Induced representations of
reductive $p$-adic groups I, {\it Ann. Sci. Ec. Norm. Super.} 10 (1977),
441-472.}

\refe{[C]}{L. Clozel, Characters of non-connected, reductive $p$-adic groups,
{\it Canad. J. Math.} 39 (1987), 149-167.}


\refe{[Fsym]}{Y. Flicker, On the symmetric-square. Applications of a trace
formula; {\it Trans. AMS} 330 (1992), 125-152;  Total global comparison; {\it
J. Funct. Anal.} 122 (1994), 255-278; Unit elements; {\it Pacific J. Math.}
175 (1996), 507-526; {\it On the Symmetric Square Lifting}, part I of
{\it Automorphic Representations of Low Rank Groups}, research monograph.}

\refe{[F]}{Y. Flicker, {\it Matching of orbital integrals on} GL(4) {\it
and} GSp(2), {\it Memoirs AMS} 137 (1999), 1-114; Automorphic forms on SO(4);
{\it Proc. Japan Acad.} 80 (2004), 100-104; Automorphic forms on PGSp(2);
{\it Elect. Res. Announc. AMS} 10 (2004), 39-50. http://www.ams.org/era/;
{\it Lifting Automorphic Representations of} PGSp(2) {\it and} SO(4) {\it to}
PGL(4), part I of {\it Automorphic Forms and Shimura Varieties of} PGSp(2),
World Scientific, 2005.}

\refe{[FK]}{Y. Flicker, D. Kazhdan, On the symmetric-square. Unstable
local transfer, {\it Invent. Math.} 91 (1988), 493-504.}

\refe{[FZ3]}{Y. Flicker, D. Zinoviev, On the symmetric-square. Unstable
twisted characters, {\it Israel J. Math.} 134 (2003), 307-315.}

\refe{[FZ4]}{Y. Flicker, D. Zinoviev, Twisted character of a small
representation of PGL(4), {\it Moscow Math. J.} 4 (2004), 333-368.}

\refe{[H]}{Harish-Chandra, {\it Admissible invariant distributions on
reductive} $p$-{\it adic groups}, notes by S. DeBacker and P. Sally, AMS Univ.
Lecture Series 16 (1999); see also: {\it Queen's Papers in Pure and Appl.
Math.} 48 (1978), 281-346.}

\refe{[K]}{D. Kazhdan, On liftings, in {\it Lie Groups Representations II},
Springer Lecture Notes on Mathematics 1041 (1984), 209-249.}

\refe{[KS]}{R. Kottwitz, D. Shelstad, {\it Foundations of Twisted Endoscopy},
Asterisque 255 (1999), vi+190 pp.}

\refe{[LN]}{R. Lidl, H. Niederreiter, {\it Finite Fields}, Cambridge Univ.
Press, 1997.}

\refe{[Wa]}{J.-L. Waldspurger, Sur les int\'egrales orbitales tordues pour
les groupes lin\'eaires: un lemme fondamental, {\it Canad. J. Math.} 43
(1991), 852-896.}

\refe{[W]}{R. Weissauer, A special case of the fundamental lemma, preprint.}



\enddocument